# The football {5, 6, 6} and its geometries: from a sport tool to fullerens and further

**Emil Molnár**[1]
Budapest University of Technology and Economics,
Institute of Mathematics, Department of Geometry
H-1521 Budapest, Hungary
emolnar@math.bme.hu

**István Prok**
Budapest University of Technology and Economics,
Institute of Mathematics, Department of Geometry
H-1521 Budapest, Hungary
prok@math.bme.hu

**Jenő Szirmai**
Budapest University of Technology and Economics,
Institute of Mathematics, Department of Geometry
H-1521 Budapest, Hungary
szirmai@math.bme.hu

## Abstract

This presentation starts with the regular polygons, of course, then with the Platonic and Archimedean solids. The latter ones are whose symmetry groups are transitive on the vertices, and in addition, whose faces are regular polygons (see only I. Prok's home page [11] for them). Then there come these symmetry groups themselves (starting with the cube and octahedron, of course, then icosahedron and dodecahedron). Then come the space filling properties: Namely the cube is a space filler for the Euclidean space $\mathbf{E}^3$. Then we jump for the other regular solids that cannot fil $\mathbf{E}^3$, but can hyperbolic space $\mathbf{H}^3$, a new space. These can be understood better if we start regular polygons, of course, that cannot fil $\mathbf{E}^2$ in general, but can fil the new plane $\mathbf{H}^2$, as hyperbolic or Bolyai-Lobachevsky plane. Now it raises the problem, whether a football polyhedron - with its congruent copies - fil a space. It turns out that $\mathbf{E}^3$ is excluded (it remains an open problem – for you, of course, in other aspects), but $\mathbf{H}^3$ can be filled as a schematic construction show this (Fig. 5), far from elementary. Then we mention some stories on Buckminster Fuller, an architect, who imagined first time fullerens as such crystal structures. Many problems remain open, of course, we are just in the middle of living science.

---

[1] Corresponding author

## 1. Introduction

After the above abstract we start with the regular *n*-sided polygons and its symmetry groups, generated by two line reflections **a** and **b** with relations $\mathbf{1} = \mathbf{a}^2 = \mathbf{b}^2 = (\mathbf{ab})^n$, i.e. the edge bisector line *a* and angle bisector *b* intersect each other in the centre *O* of the *n*-gon *in angle* $\pi/n$. Then come the 5 regular Platonic solids with symbols (*p*,*q*): (3,3) = tetrahedron (self dual), (3,4) = octahedron, (4,3) = cube, (3,5) = icosahedron, (5,3) = dodecahedron.

These can also be derived by an elementary observation, which asserts that the angle of a regular Euclidean *p*-gon is $(p\text{-}2)\pi/p$ (the angle sum of a (say convex) *p*-gon is $(p\text{-}2)\pi$, since it can be divided by *p*-3 diagonals, from any vertex, into *p*-2 triangles, each having an angle sum $\pi = 180°$). So *q* pieces of them meet at a vertex (with non-plane vertex figure), if

$$\frac{q(p-2)\pi}{p} < 2\pi \quad \Leftrightarrow \text{ (equivalent to) } \frac{1}{2} < \frac{1}{p} + \frac{1}{q},$$

a necessary condition. This will be sufficient if we think of a spherical triangle (so-called *characteristic triangle*) with angles $\pi/2$, $\pi/p$ and $\pi/q$ (with angle sum larger than $\pi$), generating the above vertex figure.

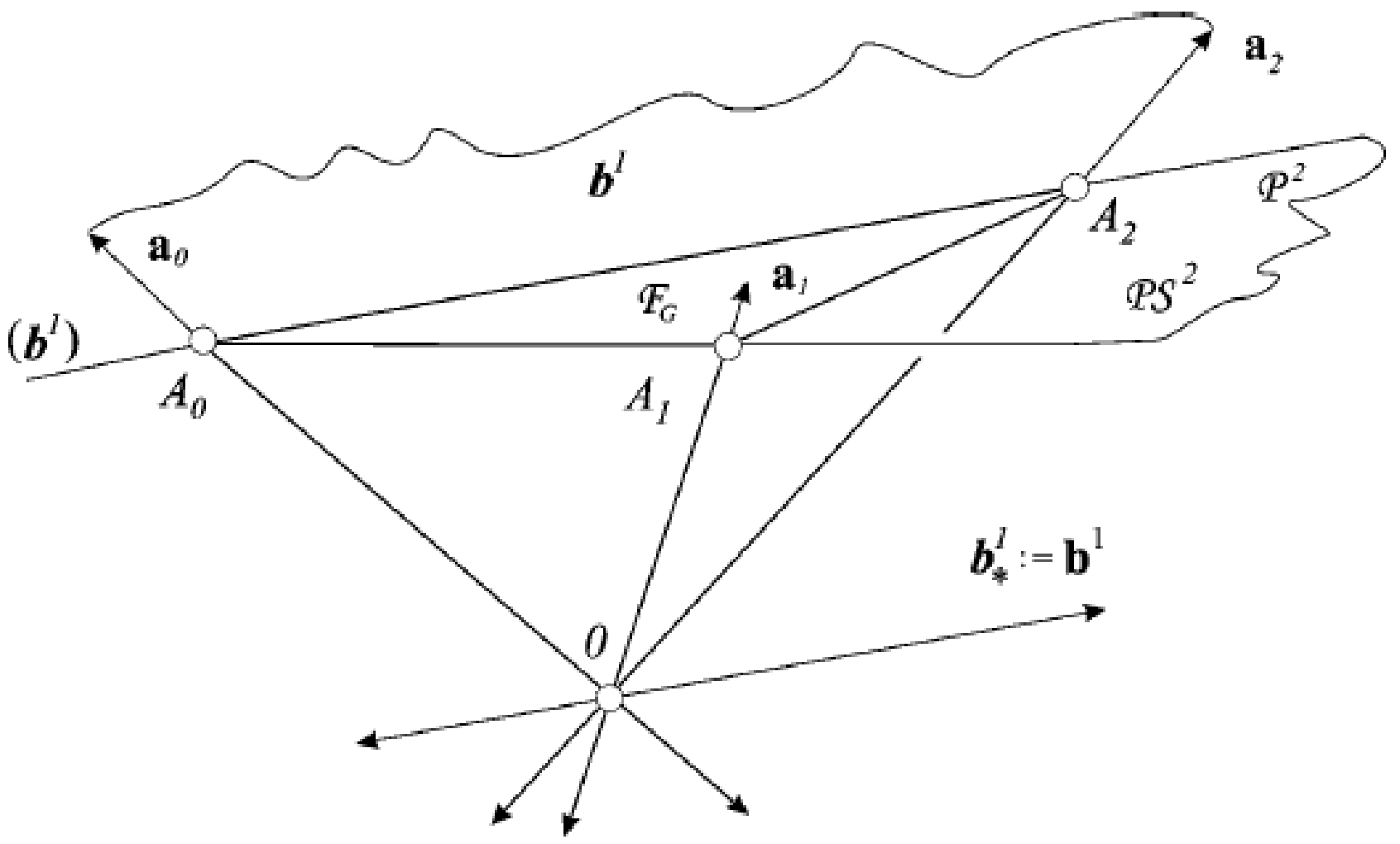

Figure 1: A projective coordinate triangle to our fundamental triangle $A_0\, A_1\, A_2$. The plane by form class $\boldsymbol{b}^1$ (e.g.) describes the line $A_0A_2$. The point by polar vector $\boldsymbol{b}^1{}_* = \mathbf{b}^1 = b^{10}\mathbf{a}_0 + b^{11}\mathbf{a}_1 + b^{12}\mathbf{a}_2 =: b^{1j}\mathbf{a}_j$ is its pole.

In more general extension, (*p*,*q*) and (*q*,*p*) are dual pairs: as the octahedron and cube, icosahedron and dodecahedron, above. To the symmetry group **G** we have 3 generating plane reflections, denoted by $\boldsymbol{b}^0$, $\boldsymbol{b}^1$, $\boldsymbol{b}^2$ (with the corresponding planes $b^i$, $i \in \{0, 1, 2\}$), and defining relations: $\mathbf{1} = \boldsymbol{b}^0\boldsymbol{b}^0 = \boldsymbol{b}^1\boldsymbol{b}^1 = \boldsymbol{b}^2\boldsymbol{b}^2 = (\boldsymbol{b}^0\boldsymbol{b}^1)^p = (\boldsymbol{b}^0\boldsymbol{b}^2)^2 = (\boldsymbol{b}^1\boldsymbol{b}^2)^q$, as the Coxeter-Schläfli diagram

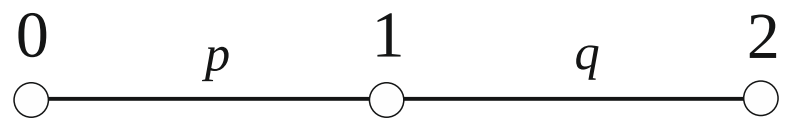

(1.1)

indicates this with nodes o for the reflection planes, and the *branches* $p$ and $q$ for the relations above. If two nodes are not connected then their reflections commute (with their orthogonal planes, the product order is 2 then). Thus the tetrahedron group is of order 24, as the symmetry group **G** of the platonic solid are in general of order $|\mathbf{G}| = 4\pi /(\pi/p + \pi/q - \pi/2)$. Namely, the surface area $4\pi/R^2$ of the sphere $\mathbf{S}^2$ of radius $R$ is divided by the area (with angle excess $R^2[(\pi/p + \pi/q + \pi/2) - \pi]$) of the fundamental triangle of **G** (see the above elementary observations as well). For octahedron and cube we have 48 symmetry elements, for icosahedron and dodecahedron we have 120 elements for the symmetry group **G**.

At the same time we can introduce the so-called Coxeter-Schläfli matrix

$$(b^{ij}) = \langle \boldsymbol{b}^i, \boldsymbol{b}^j \rangle = \begin{bmatrix} 1 & -\cos\left(\frac{\pi}{p}\right) & 0 \\ -\cos\left(\frac{\pi}{p}\right) & 1 & -\cos\left(\frac{\pi}{q}\right) \\ 0 & -\cos\left(\frac{\pi}{q}\right) & 1 \end{bmatrix} = (\cos(\pi - \beta^{ij})). \qquad (1.2)$$

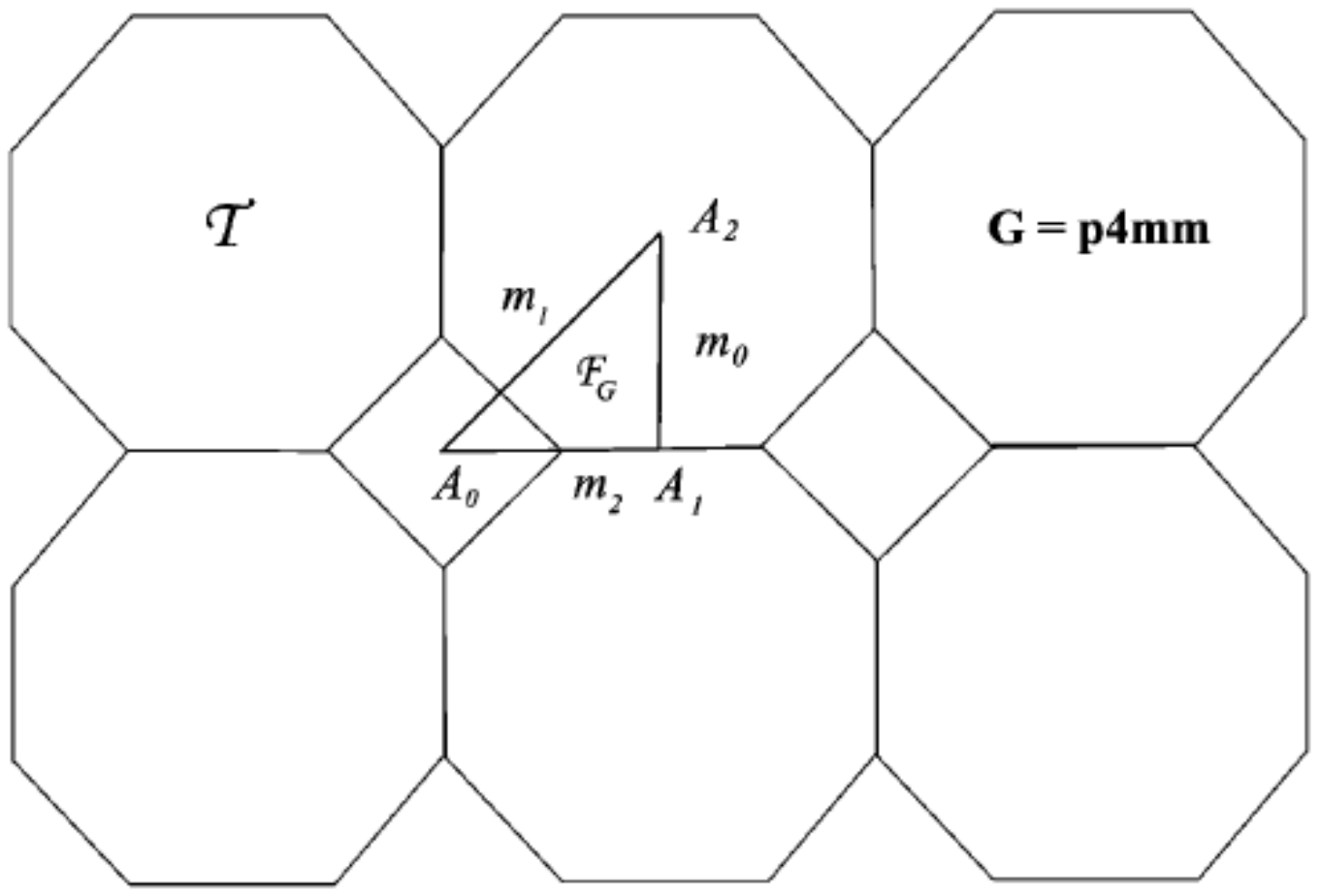

Figure 2: A well-known Archimedean tiling in the Euclidean plane with its fundamental triangle.

This is derived from the formal scalar products of the linear forms (normal unit vectors) $\boldsymbol{b}^i$, ordered to the side lines $b^i$ of the characteristic triangle $A_0A_1A_2$ with angles $\beta^{01} = \pi/p$, $\beta^{02} = \pi/2$, $\beta^{12} = \pi/q$. By convention, any side line $b^i$ above closes angle $\pi$ with itself, it lies opposite to vertex $A_i$.

Fig. 1 shows symbolically a (2+1)–dimensional picture to the projective triangle (simplex) coordinate system, also for the later higher dimensional analogue in a ($d$+1)–dimensional vector space $\mathbf{V}^{d+1}$ and its dual $\boldsymbol{V}_{d+1}$. Here $O$ denotes the origin from where vectors $\mathbf{a}_i = OA_i$ ($i \in \{0, 1, 2\}$ point to the vertices of triangle $A_0A_1A_2$, forms as normal vectors $\boldsymbol{b}^j$ are placed *to the side plane* $b^j$. See also the later Fig. 2 for an *Archimedean tiling* in $\mathbf{E}^2$. Thus

$$\mathbf{a}_i\boldsymbol{b}^j = \delta_i^j$$

(the *Kronecker symbol*, $i, j \in \{0, 1, 2\}$) indicates the incidence relations. By the way, the inverse matrix $(a_{ij}) = (b^{ij})^{-1}$ of the above Coxeter-Schläfli matrix just serves the distance metrics of the triangle $A_0A_1A_2$. The side length $A_iA_j$ can be expressed from spherical angle by

$$\cos\left(\frac{A_iA_j}{R}\right) = \frac{a_{ij}}{\sqrt{a_{ii}a_{jj}}}. \tag{1.3}$$

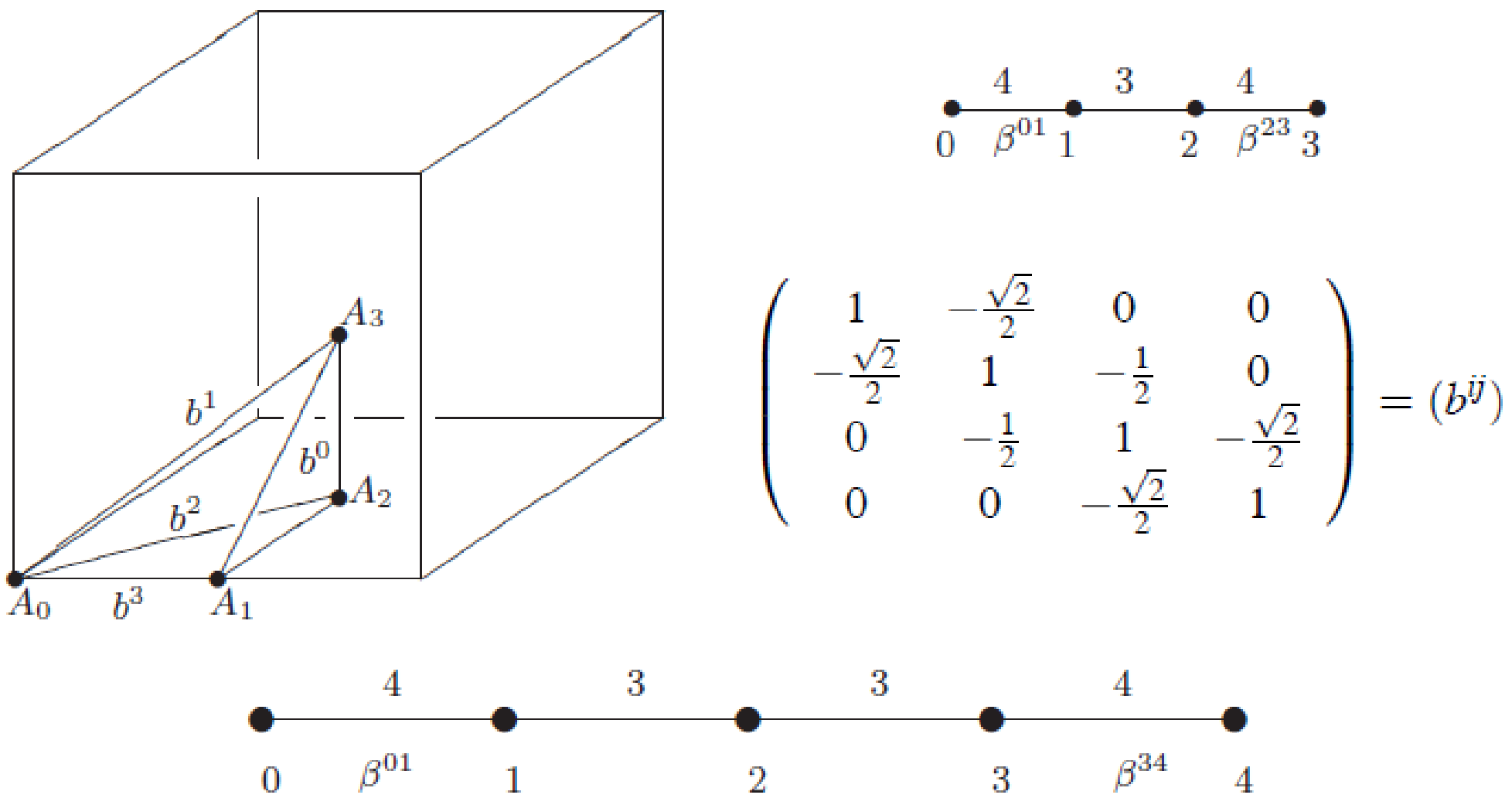

Figure 3: Cube tiling in $\mathbf{E}^3$ and symbols for it. Coxeter-Schläfli diagram for the $\mathbf{E}^4$ cube tiling.

This is related now to the spherical geometry of the Platonic solid above, where

$$\det(b^{ij}) = 1 - \cos^2(\pi/p) - \cos^2(\pi/q) > 0. \tag{1.4}$$

E.g. for the cube $(p, q) = (4, 3)$ we get $\det(b^{ij}) = 1 - \frac{1}{2} - \frac{1}{4} = \frac{1}{4}$. This is more critical for the dodecahedron (or icosahedron): $\det(b^{ij}) = 1 - \cos^2(\pi/5) - \frac{1}{4} = \frac{3-\sqrt{5}}{8} \approx 0.09549 > 0$.

## 2. Euclidean and non-Euclidean mosaics

If above $\det(b^{ij}) = 0$ in formula (1.4), i.e. $(p, q) = (3, 6), (6, 3), (4, 4)$, then we get mosaic or tiling in Euclidean plane $\mathbf{E}^2$. The latter one is the well-known squared paper. These follow also by the above elementary observations.

But we can assume in formula (1.4) also $\det(b^{ij}) < 0$, i.e. then we get an infinite series for $(p, q) = (3, 7), (7, 3), (3, 8), (8, 3), \ldots, (4, 5), (5, 4), (4, 6), (6, 4), \ldots$ . Above, we can imagine characteristic triangles with angles $\pi/p$, $\pi/q$, $\pi/2$, their sum is smaller than $\pi = 180°$ as it holds for Euclidean plane $\mathbf{E}^2$. We have obtained a new geometry, the so-called hyperbolic or Bolyai-Lobachevsky plane, denoted by $\mathbf{H}^2$. This geometry was discovered and elaborated first, approximately in the same time and independently, by the Hungarian János Bolyai and the Russian Nikolai Ivanović Lobačevskii in 1820's years.

The area of a triangle with angles $\alpha$, $\beta$, $\gamma$ is equal (proportional) to the defect $\pi - (\alpha + \beta + \gamma)$ in $\mathbf{H}^2$. (This was observed also by the fore-runners of the new geometry, so by Carl Friedrich Gauss, who played so important role in the life of Farkas (Wolgang) and János Bolyai, see e.g. [9].) Thus the starting tilings (3, 7), (7, 3) have just a minimal characteristic triangle with area $\pi/2 - \pi/3 - \pi/7 = \pi/42$.

It turned out, that the formulas of spherical geometry $\mathbf{S}^2$ become true formulas in $\mathbf{H}^2$ if we substitute imaginary radius $ki = R$ into the spherical formulas, e.g. into (1.1) above ($i = \sqrt{-1}$ is the imaginary complex unit as usual). Thus, so-called hyperbolic functions come into the play, etc. So we get a unifying concept, *absolute geometry* in the sense of János Bolyai, a joint kernel of geometries $\mathbf{E}^2$, $\mathbf{S}^2$, $\mathbf{H}^2$, and later in any dimension $d$ for $\mathbf{E}^d$, $\mathbf{S}^d$, $\mathbf{H}^d$, where the analogy does not remain so easy, and we have many open problems.

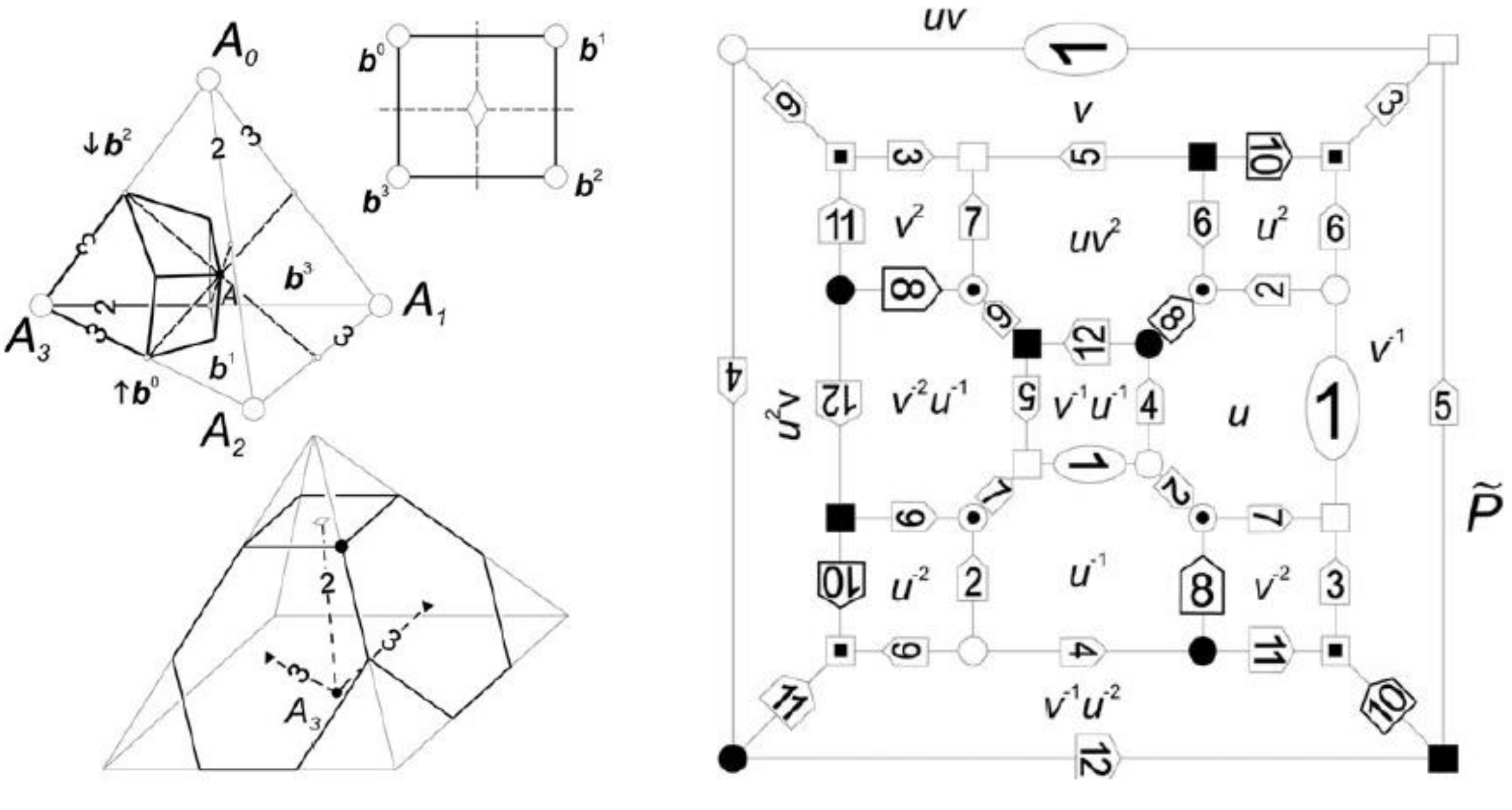

Figure 4: The Archimedean solid, {4, 6, 6} as truncated octahedron and its Euclidean manifold $\mathbf{E}^3/\mathbf{P}2_12_12_1 = \tilde{P}$.

However, for the above regular tilings ($p$, $q$) we have a unified theory (see e.g. [2] and our Sect.5), where reflection groups (so-called Coxeter groups, as above in formulas (1.1-4)), play important roles.

An interesting topic is the so-called *Archimedean or vertex transitive tilings* (mosaics, see Fig. 2) by polygons, where the symmetry group **G** of the tiling $T$ acts *transitively* on the vertices (i.e. all vertices form one equivalence class under the symmetry group $\mathbf{G}(T)$. In addition, for nice pictures, we assume that the tiles of $T$ are regular polygons (they are possibly non-congruent, of course).

The above series of Coxeter reflection groups, denoted by $[2, p, q]$ with fundamental triangle of angles $\pi/2$, $\pi/p$, $\pi/q$, each provides nice Archimedean tiling. E.g. {4, 6, 6}, where square-hexagon-hexagon meet at every vertex, is the truncated octahedron as a space filler polyhedron of Euclidean space $\mathbf{E}^3$. Analogously {4, 8, 8}, where square-octagon-octagon meet at every vertex, is a very popular pavement of our Euclidean streets (Fig. 2, 4). Moreover, {4, 10, 10} will be a hyperbolic tiling in $\mathbf{H}^2$, see also [4].

Special interest deserves the football polyhedron {5, 6, 6}, where pentagon-hexagon-hexagon meet at every vertex (Fig. 5). Thus, we obtain 12 pentagons and 20 hexagons with 60 vertices at all. This polyhedron can be obtained from the regular icosahedron (with group [2, 3, 5] above) by truncating its 12 vertices each derives a regular pentagon, so that the 20 triangles become 20 hexagons. This football "plays important role in our life", of course, and it turned out, this can play new roles in the structure of new materials as fullerens in crystallography, may be in non-Euclidean crystallography as follows in the next sections.

## 3. Space filler polyhedra in Euclidean and non-Euclidean spaces

We know that the cube is a space filler polyhedron, i.e. we can fil Euclidean space $\mathbf{E}^3$ with its congruent copies, face-to-face without gaps and overlaps. Namely, this can be derived by the arguments close to that of the introduction at the Coxeter–Schläfli diagram and matrix, now especially with $(p, q, r) = (4, 3, 4)$ with some extension:

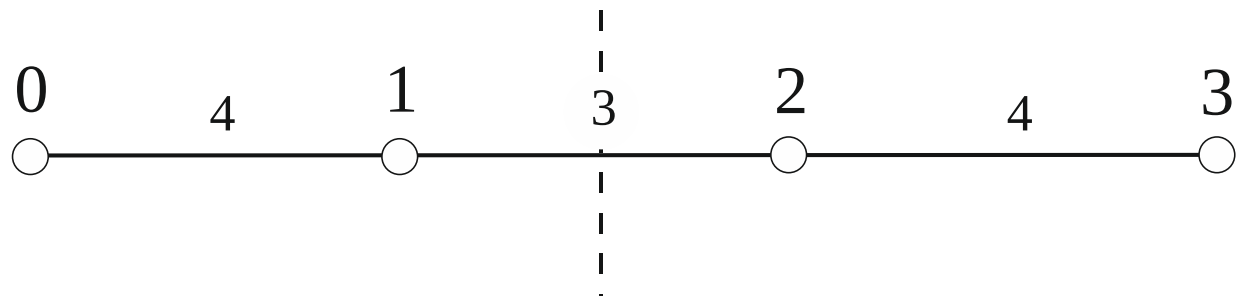

(3.1)

$$(b^{ij}) = \langle \boldsymbol{b}^i, \boldsymbol{b}^j \rangle = \begin{bmatrix} 1 & -\cos\left(\frac{\pi}{p}\right) & 0 & 0 \\ -\cos\left(\frac{\pi}{p}\right) & 1 & -\cos\left(\frac{\pi}{q}\right) & 0 \\ 0 & -\cos\left(\frac{\pi}{q}\right) & 1 & -\cos\left(\frac{\pi}{r}\right) \\ 0 & 0 & -\cos\left(\frac{\pi}{r}\right) & 1 \end{bmatrix} = (\cos(\pi - \beta^{ij})). \quad (3.2)$$

That means (Fig. 3), we have a characteristic tetrahedron (simplex) $A_0\, A_1\, A_2\, A_3 = b^0\, b^1\, b^2\, b^3$ with the cube centre $A_3$, with a face centre $A_2$, with midpoint $A_1$ of an edge incident to the previous face, and a vertex $A_0$ of the previous edge. Furthermore, $b^i = A_j\, A_k\, A_l$ with $\{i, j, k, l\} = \{0, 1, 2, 3\}$ hold. We apply real (left) vector space $\mathbf{V}^4$ for points and its (right) dual space $\boldsymbol{V}_4$ for planes up to projective equivalences ~ as usual. Vector coefficients are written in rows on the left to column basis; dual forms with row basis have column coefficients on the right. Forms act on the right on vectors.

Now the characteristic simplex has a half-turn symmetry $0 \leftrightarrow 3$, $1 \leftrightarrow 2$ as indicated in the extended diagram (3.1). The simple reflection diagram describes the crystallographic space group 221. **Pm3m** with so-called primitive cube lattice [3]. The (with half-turn) extended diagram describes the space group 229. **Im3m** with body centred (*innenzentriert* in German) cube lattice. The last assertions indicate that our topic has important applications in crystallography, describing atomic or molecular structures, etc. You find [3] in the Internet for free download. We will extend these concepts to non-Euclidean geometry as well.

## 4. A Euclidean space form, as typical example

For the later generalization, we analyse Fig. 4 with the truncated octahedron {4, 6, 6} as Archimedean solid that can fil Euclidean space $\mathbf{E}^4$ with its congruent copies. To this a famous space filler tetrahedron (simplex) $A_0\, A_1 A_2\, A_3 = b^0\, b^1\, b^2\, b^3$, the so-called *sphenoid*, with extended diagram in Fig.4 (to the space group 224. **Pn3m)** and Coxeter–Schläfli matrix

$$(b^{ij}) = (\langle \boldsymbol{b}^i, \boldsymbol{b}^j \rangle) = \begin{bmatrix} 1 & -\frac{1}{2} & 0 & -\frac{1}{2} \\ -\frac{1}{2} & 1 & -\frac{1}{2} & 0 \\ 0 & -\frac{1}{2} & 1 & -\frac{1}{2} \\ -\frac{1}{2} & 0 & -\frac{1}{2} & 1 \end{bmatrix} = (\cos(\pi - \beta^{ij})) \qquad (4.1)$$

play important roles. That means, the simplex has two opposite rectangles $\beta^{02} = \beta^{13} = \pi/2$ as dihedral angles, the remaining four dihedral angles are $\pi/3$. Moreover, the simplex has 3 half-turns, their axes connect the corresponding opposite edge midpoints and meet at the simplex centre $A$. The images of $A$ will just form the Archimedean vertex class equivalent under the 24 regular tetrahedron symmetries around $A_3$. Thus, we get the truncated octahedron {4, 6, 6} denoted by $P$ and its so-called Schlegel diagram in Fig. 4 right in $\mathbf{E}^2$.

Now we organize new face pairing isometries and new space filler tiling with $\tilde{P}$ = {4, 6, 6} as fundamental polyhedron under a new fixed-point-free space group (finally it will be 19. $\mathbf{P2_12_12_1}$ by [3]).

Our method will be very general, and we apply it in the following.

Look at the sphenoid in Fig. 4 left that 3 truncated octahedron meet at any edge of the {4, 6, 6} tiling, along square-hexagon, hexagon-hexagon, hexagon-square. Thus in Fig. 4 right, we want to construct a new group with fixed-point-free action, where any point will have a ball-like neighbourhood without its other image point in that ball. We start a directed edge class consisting of three edges, denoted say by 1. The first 1 is defined on the boundary of a square and a hexagon denoted by $u^{-1}$, then the second 1 is on the **u**-image hexagon $u$, while hexagon $v^{-1}$ follows on its other side, then the third 1 is placed on the **v**-image hexagon $v$ with a square on its other side. Now the first square will be denoted by $(uv)^{-1} = v^{-1}u^{-1}$, the last square (now with the outer „infinite" face) by $uv$ because the side pairing mapping with usual conventions:

$$\mathbf{uv}:\ (uv)^{-1} = v^{-1}u^{-1} \to\ uv \ \text{ and } (\mathbf{uv})^{-1} = \mathbf{v}^{-1}\mathbf{u}^{-1}:\ uv \to (uv)^{-1} = v^{-1}u^{-1}. \qquad (4.2)$$

Imagine our new tiling as representing the elements of our new group **G** (= $\mathbf{P2_12_12_1}$) by the images of a starting fundamental domain $\tilde{P}^{\mathbf{1}}$ (the identity domain). The **u** -image domain $\tilde{P}^{\mathbf{u}}$ lies just besides face $u$ with its face $(u^{-1})^{\mathbf{u}}$. Similarly, the $\mathbf{u}^{-1}$-image domain $\tilde{P}^{\mathbf{u}^{-1}}$ lies besides face $u^{-1}$ with face $(u)^{\mathbf{u}^{-1}}$. Thus, to the first edge 1 and face $u^{-1}$ of $\tilde{P}$, we find three edge-domains each between two faces:

$$(v^{-1}u^{-1})|\,\tilde{P}|(u^{-1}) \ \text{ then } \ (u)^{\mathbf{u}^{-1}}|(\tilde{P})^{\mathbf{u}^{-1}}|(v^{-1})^{\mathbf{u}^{-1}} \ \text{ then } \ (v)^{\mathbf{v}^{-1}\mathbf{u}^{-1}}|(\tilde{P})^{\mathbf{v}^{-1}\mathbf{u}^{-1}}|(uv)^{\mathbf{v}^{-1}\mathbf{u}^{-1}} \qquad (4.3)$$

then cyclically comes the first identity edge domain, now with a formal general rule, as **uv**-image of $(\tilde{P})^{\mathbf{v}^{-1}\mathbf{u}^{-1}}$, i.e. $[(\tilde{P})^{\mathbf{v}^{-1}\mathbf{u}^{-1}}]^{\mathbf{uv}} = \tilde{P}$.

Consider the edge between $u$ and $u^{-1}$, numbered by 2 as an arrow. The mappings **u** and $\mathbf{u}^{-1}$ fixed before, carry this edge 2 into two other edges denoted also by 2: first to the edge between $u$ and a square face, denoted by $u^2$, second to that edge between $u^{-1}$and the square $u^{-2}$. So, the face pairing isometries to the edge class 2:

$$\mathbf{u}^{2}:\ u^{-2} \to u^{2} \ \text{ and } \ \mathbf{u}^{-2}:\ u^{2} \to u^{-2} \qquad (4.4)$$

have also been introduced.

Then we obtain a straightforward procedure with 12 edge classes: either we get a new face pairing isometry expressed by generators **u** and **v**, mapping $\tilde{P}$ onto an adjacent image; or we get a trivial relation, going around an edge in the tiling, e.g. at 6: $\mathbf{u}^2\mathbf{v}\mathbf{v}^{-1}\mathbf{u}^{-2} = \mathbf{1}$, as consequence; or we get a non-trivial so-called defining relation, e.g. at 8: $\mathbf{v}^2\mathbf{u}\mathbf{v}^2\mathbf{u}^{-1} = \mathbf{1}$; or we get a consequence of the former defining relations.

Now it turns out, as a lucky situation, that the procedure goes smoothly to the end, without contradiction. Every edge class has 3 edges with the same rules of adjacencies, and two generators are sufficient. So we obtain a presentation of our space group named 19. $\mathbf{P2_12_12_1}$ in [3] (see also author's further papers in references, especially [5], [7], [10]):

$$\mathbf{G} = \mathbf{P2_12_12_1} = (\mathbf{u}, \mathbf{v} \mid \mathbf{v}^2\mathbf{u}\mathbf{v}^2\mathbf{u}^{-1} = \mathbf{u}^2\mathbf{v}^{-1}\mathbf{u}^2\mathbf{v} = \mathbf{1}). \qquad (4.5)$$

The vertices of $\tilde{P}$ are divided into 6 equivalence classes, four vertices in each class. To every vertex classes join 4 edge classes (as chemical bounds between atoms) as the next formula shows:

$$\circ\ (1, 2, 4, 6); \quad \square\ (1, 3, 5, 7); \quad \bullet\ (4, 8, 11, 12);$$
$$\blacksquare\ (5, 9, 10, 12); \quad \odot\ (2; 7; 8; 9); \quad \boxdot\ (3; 6; 10; 11) \qquad (4.6)$$

*Table 1. The truncated octahedron, as fundamental domain for the following Euclidean space groups [3], with its face pairings given up to its symmetries (computer classification by István Prok and Zsanett Szuda (BME Math MSc student)).*

| No. | Name | # |
|---|---|---|
| 1 | **P1** | 1 |
| 2 | $\mathbf{P\bar{1}}$ | 2 |
| 4 | $\mathbf{P2_1}$ | 2 |
| 5 | **C2** | 3 |
| 7 | **Pc** | 2 |
| 9 | **Cc** | 3 |
| 13 | **P2/c** | 3 |
| 14 | $\mathbf{P2_1/c}$ | 10 |
| 15 | **C2/c** | 8 |
| 16 | **P222** | 1 |
| 18 | $\mathbf{P2_12_12}$ | 2 |
| 19 | $\mathbf{P2_12_12_1}$ | 2 |
| 20 | $\mathbf{C222_1}$ | 2 |
| 23 | **I222** | 2 |
| 27 | **Pcc2** | 2 |
| 29 | $\mathbf{Pca2_1}$ | 3 |
| 32 | **Pba2** | 1 |
| 33 | $\mathbf{Pna2_1}$ | 4 |

| No. | Name | # |
|---|---|---|
| 34 | **Pnn2** | 1 |
| 41 | **Aea2** | 3 |
| 43 | **Fdd2** | 1 |
| 45 | **Iba2** | 3 |
| 50 | **Pban** | 3 |
| 54 | **Pcca** | 3 |
| 56 | **Pccn** | 3 |
| 60 | **Pbcn** | 6 |
| 61 | **Pbca** | 3 |
| 68 | **Ccce** | 2 |
| 70 | **Fddd** | 1 |
| 76 | $\mathbf{P4_1}$ | 1 |
| 78 | $\mathbf{P4_3}$ | |
| 77 | $\mathbf{P4_2}$ | 1 |
| 81 | $\mathbf{P\bar{4}}$ | 1 |
| 82 | $\mathbf{I\bar{4}}$ | 1 |
| 86 | $\mathbf{P4_2/n}$ | 1 |
| 88 | $\mathbf{I4_1/a}$ | 2 |

| No. | Name | # |
|---|---|---|
| 91 | $\mathbf{P4_122}$ | 1 |
| 95 | $\mathbf{P4_322}$ | |
| 92 | $\mathbf{P4_12_12}$ | 2 |
| 96 | $\mathbf{P4_32_12}$ | |
| 94 | $\mathbf{P4_22_12}$ | 1 |
| 98 | $\mathbf{I4_122}$ | 1 |
| 106 | $\mathbf{P4_2bc}$ | 1 |
| 110 | $\mathbf{I4_1cd}$ | 1 |
| 112 | $\mathbf{P\bar{4}2c}$ | 1 |
| 114 | $\mathbf{P\bar{4}2_1c}$ | 2 |
| 117 | $\mathbf{P\bar{4}b2}$ | 1 |
| 120 | $\mathbf{I\bar{4}c2}$ | 1 |
| 122 | $\mathbf{I\bar{4}2d}$ | 1 |
| 144 | $\mathbf{P3_1}$ | 1 |
| 145 | $\mathbf{P3_2}$ | |
| 152 | $\mathbf{P3_121}$ | 1 |
| 154 | $\mathbf{P3_221}$ | |

With the scalar product $\langle \boldsymbol{b}^i, \boldsymbol{b}^j \rangle = b^{ij}$ of matrix (4.1) we can formally define a quadratic form (with Einstein-Schouten index convention for summing):

$$w_i b^{ij} w_j = w_0w_0 - w_0w_1 - w_0w_3 + w_1w_1 - w_1w_2 + w_2w_2 - w_2w_3 + w_3w_3 =$$
$$= (w_0 - 1/2\, w_1 - 1/2\, w_3)^2 + 3/4\,(w_1 - 2/3\, w_2 - 1/3\, w_3)^2 + 2/3\,(w_2 - w_3)^2 \qquad (4.7)$$

with sum of three positive squares. We say that the scalar product of matrix (4.1) is of signature (+,+,+,0), and this is characteristic for Euclidean geometry $\mathbf{E}^3$. The determinant of matrix (4.1) equals to 0 (zero). But take $(p, q, r) = (5, 3, 5)$ in matrix (3.2) instead of diagram (3.1). Then the determinant will be negative and we would obtain a negative square summand besides of the three positive square summand in the quadratic form, analogue to (4.7) (not uniquely).

This will be the topic in the next main sections.

However, before that, we add Table 1, as complete computer classification of all the space groups which can have the truncated octahedron as fundamental domain. Its face pairings are given up to its symmetries.

## 5. The hyperbolic space $\mathbf{H}^3$ to the football {5, 6, 6}

We start with our Fig. 5, the plane figure of the real football (again, the so-called Schlegel diagram). Maybe surprisingly, the Coxeter–Schläfli matrix and diagram provide us with all the theoretical tools, as above. Now it will be in the hyperbolic space $\mathbf{H}^3$, modelled in the projec-tive metric space $\mathcal{P}^3(\mathbf{V}^4, \boldsymbol{V}_4, \sim, \langle,\rangle)$ with a scalar product (or polarity) of above signature (+,+,+,–). For the Coxeter-Schläfli diagram in Fig. 5 we formally introduce the matrix

$$(b^{ij}) = \begin{bmatrix} 1 & -\cos\left(\frac{\pi}{5}\right) & 0 & 0 \\ -\cos\left(\frac{\pi}{5}\right) & 1 & -\cos\left(\frac{\pi}{3}\right) & 0 \\ 0 & -\cos\left(\frac{\pi}{3}\right) & 1 & -\cos\left(\frac{\pi}{5}\right) \\ 0 & 0 & -\cos\left(\frac{\pi}{5}\right) & 1 \end{bmatrix}, \tag{5.1}$$

and fix basis $\boldsymbol{b}^i$ of (the dual form space) $\boldsymbol{V}_4$ to the side faces of simplex $A_0A_1A_2A_3 = b^0b^1b^2b^3$ with the scalar product

$$\langle,\rangle:\ \boldsymbol{V}_4 \times \boldsymbol{V}_4 \to \mathbf{R} \text{ (the real field) by } \langle \boldsymbol{b}^i, \boldsymbol{b}^j \rangle = b^{ij} \text{ with the matrix (5.1).} \tag{5.2}$$

This scalar product has the signature (+,+,+,–) now, since

$$B := \det(b^{ij}) = \left(1 - \cos^2\frac{\pi}{5}\right)^2 - \cos^2\frac{\pi}{3} < 0, \tag{5.3}$$

however, all the principal minors are positive in the chain

$$b^{00}=1,\ \begin{vmatrix} b^{00} & b^{01} \\ b^{10} & b^{11} \end{vmatrix} = 1 - \cos^2\frac{\pi}{5},\ \begin{vmatrix} b^{00} & b^{01} & b^{02} \\ b^{10} & b^{11} & b^{12} \\ b^{20} & b^{21} & b^{22} \end{vmatrix} = 1 - \cos^2\frac{\pi}{5} - \cos^2\frac{\pi}{3} > 0. \tag{5.4}$$

The inverse matrix $(a_{ij})$ of $(b^{ij})$, with $b^{ij} \cdot a_{jk} = \delta^i_k$, has great importance, since it induces the scalar product on $\mathbf{V^4}$ by linear extension

$$\langle , \rangle:\ \mathbf{V^4} \times \mathbf{V^4} \to \mathbf{R} \ \text{ by } \langle \mathbf{a}_i, \mathbf{a}_j \rangle = a_{ij}\,, \tag{5.5}$$

where $\{\mathbf{a}_j\}$ is just the dual basis to $\{\boldsymbol{b}^i\}$ defined by $\mathbf{a}_j \boldsymbol{b}^i = \delta^i_j$. Geometrically, the vectors $\mathbf{a}_j$ represent the vertices $A_j$ of the Lanner simplex $F_L$ whose side planes $m_i$ are described by forms $\boldsymbol{b}^i$. In general, the vectors $\mathbf{x}$ in the cone

$$C = \{\mathbf{x} \in \mathbf{V}^4 : \langle \mathbf{x}; \mathbf{x} \rangle < 0\} \tag{5.6}$$

define the proper points $(\mathbf{x})$ of the hyperbolic space $\mathbf{H^3}$ embedded in $\mathcal{P}^3$. If $(\mathbf{x})$ and $(\mathbf{y})$ are proper points, with $\langle \mathbf{x} ; \mathbf{y} \rangle < 0$, then their distance $d(\mathbf{x}, \mathbf{y})$ is defined by

$$\cosh \frac{d}{k} = \frac{-\langle \mathbf{x}; \mathbf{y} \rangle}{\sqrt{\langle \mathbf{x}; \mathbf{x} \rangle \langle \mathbf{y}; \mathbf{y} \rangle}} \geq 1, \tag{5.7}$$

where $k=\sqrt{-1/K}$ is the metric constant of $\mathbf{H^3}$ (of constant negative curvature $K$). The forms $\boldsymbol{u}$ in the complementary cone

$$C^* = \{\boldsymbol{u} \in \boldsymbol{V}_4 : \langle \boldsymbol{u}; \boldsymbol{u} \rangle > 0\} \tag{5.8}$$

define the proper planes $(\boldsymbol{u})$ of $\mathbf{H^3}$. Suppose that $(\boldsymbol{u})$ and $(\boldsymbol{v})$ are proper planes. They intersect in a proper straight line, iff $\langle \boldsymbol{u}; \boldsymbol{u} \rangle \langle \boldsymbol{v}; \boldsymbol{v} \rangle - \langle \boldsymbol{u}; \boldsymbol{v} \rangle^2 > 0$. One of their angles $\alpha(\boldsymbol{u} ; \boldsymbol{v})$ can be defined by

$$\cos \alpha = \frac{\langle \boldsymbol{u}; \boldsymbol{v} \rangle}{\sqrt{\langle \boldsymbol{u}; \boldsymbol{u} \rangle \langle \boldsymbol{v}; \boldsymbol{v} \rangle}}, \quad 0 < \alpha < \pi. \tag{5.9}$$

The scalar products introduced in $\boldsymbol{V}_4$ and $\mathbf{V}^4$, respectively, allow us to define a bijective polarity between vectors and forms:

$$\begin{aligned} &({}_*): \boldsymbol{V}_4 \to \mathbf{V}^4; \quad \boldsymbol{u} \to \boldsymbol{u}_* =: \mathbf{u} \text{ by} \\ &\quad \mathbf{u}\boldsymbol{v} := \langle \boldsymbol{u}; \boldsymbol{v} \rangle \text{ for every } \boldsymbol{v} \in \boldsymbol{V}_4, \\ &({}^*): \mathbf{V}^4 \to \boldsymbol{V}_4; \quad \mathbf{x} \to \mathbf{x}^* =: \boldsymbol{x} \text{ by} \\ &\quad \mathbf{y}\boldsymbol{x} := \langle \mathbf{y}; \mathbf{x} \rangle \text{ for every } \mathbf{y} \in \mathbf{V}^4. \end{aligned} \tag{5.10}$$

So we get a geometric polarity between points and planes: proper points have improper polars, proper planes have improper poles. A polar $(\boldsymbol{u})$ and its pole $(\mathbf{u})$ are incident iff

$$0 = \boldsymbol{u}\mathbf{u} = \langle \boldsymbol{u}; \boldsymbol{u} \rangle = \langle \mathbf{u}; \mathbf{u} \rangle \tag{5.11}$$

Such a point (**u)** is called end of $\mathbf{H}^3$ (or point on the absolute quadratic or quadric), its polar (***u*)** is a boundary plane (tangent to the absolute at the end). They form points and planes of a quadric of type (+, +, +, –).

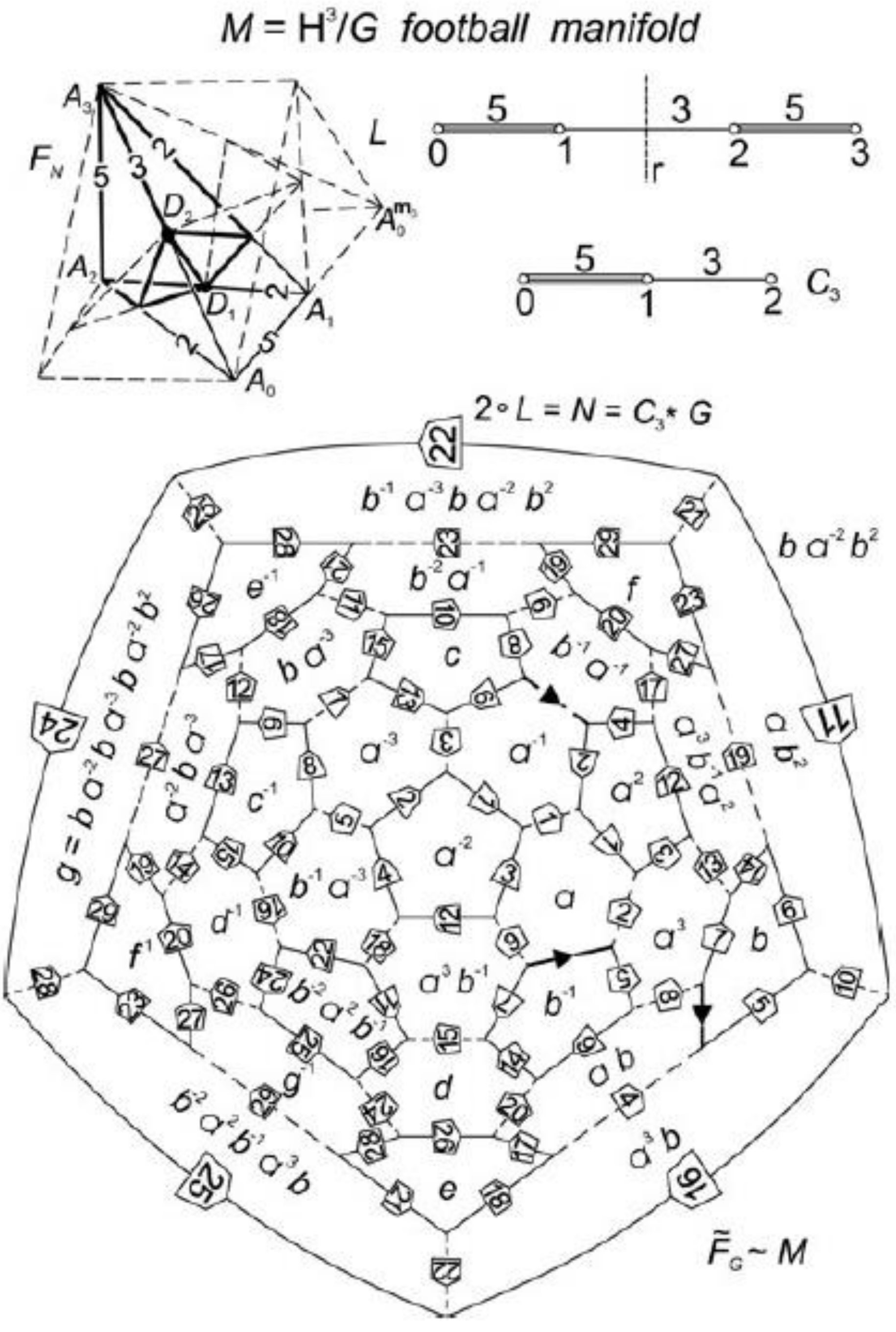

Figure 5: The hyperbolic football manifold for the Archimedean solid {5, 6, 6}.

The isometry group **G** of $\mathbf{H}^3$ can be generated by plane reflections. If (***u***) is a proper plane and (**u**) is its pole, then the reflection formulas for points (vectors) and planes (forms) are

$$x \to x - \frac{2(x\,\boldsymbol{u})}{\langle \boldsymbol{u};\boldsymbol{u}\rangle}\cdot \mathbf{u};\quad v \to v - \boldsymbol{u}\cdot\frac{2\langle v;\boldsymbol{u}\rangle}{\langle \boldsymbol{u};\boldsymbol{u}\rangle}. \tag{5.12}$$

Now, we can define the Lanner group, denoted by **L**, as a group generated by reflections in the side planes of the simplex $F_{\mathbf{L}}=A_0A_1A_2A_3$, given by the forms ($\boldsymbol{b}^j$) and their poles. As the matrix ($b^{ij}$) in (5.1) shows, the face angle between $m_0$ and $m_1$ is equal to $\frac{\pi}{5} = \beta^{01}$. Analogously hold $\beta^{02} = \frac{\pi}{2} = \beta^{03}$; $\beta^{12} = \frac{\pi}{3}$; $\beta^{13} = \frac{\pi}{2}$; $\beta^{23} = \frac{\pi}{5}$ as it has just been prescribed by the Coxeter–Schläfli diagram of **L** in Fig. 5. We can check that the subgroup generated by the plane reflections $m_0$ and $m_1$, for instance, is of order 2 × 5=10. Moreover, the subgroup $\mathbf{C}_3$, stabilizing the point $A_3$ in **L**, is just the dodecahedron group of order 120, generated by the reflections $\mathbf{m_0}$, $\mathbf{m}_1$, $\mathbf{m}_2$ in $\mathbf{H}^3$. A fundamental polyhedron of the group **L** is the simplex $A_0A_1A_2A_3=F_{\mathbf{L}}$ itself. Its points (**x**) can be characterized by

$$F_L = \left\{(\mathbf{x}) : \mathbf{x}\boldsymbol{b}^j \geq 0 \text{ for each } \boldsymbol{b}^j,\ j = 0,1,2,3\right\} \tag{5.13}$$

This is uniquely determined by the fixed point set of **L**, i. e. any other one is an **L**-image of $F_{\mathbf{L}}$. From the projective model of $\mathbf{H^3}$ we can turn to other models, e. g. to the usual Cayley-Klein model in a Euclidean ball. Our figures are shown in this model (without indicating the absolute ball).

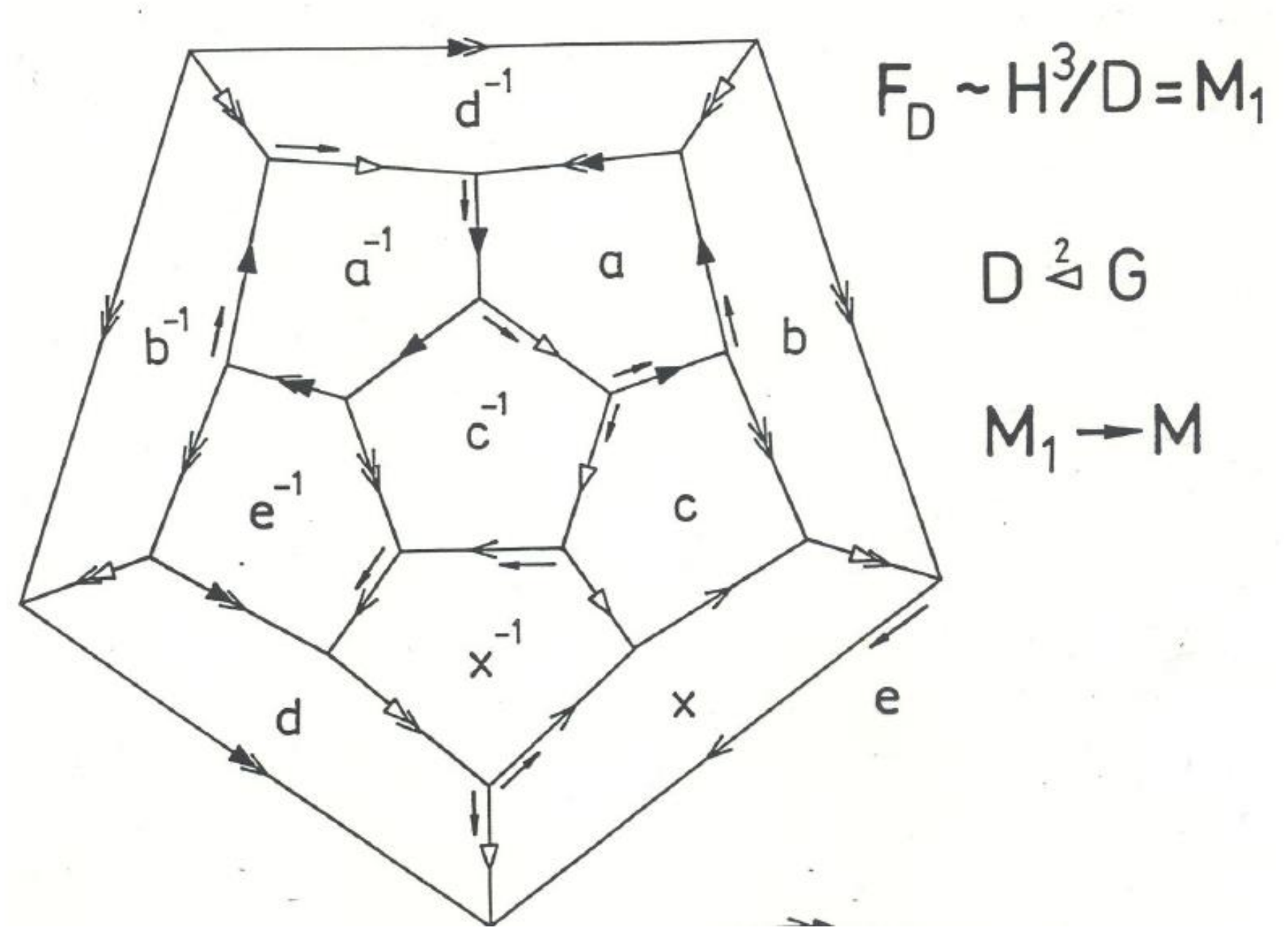

Figure 6: Twofold covering the football manifold as a hyperbolic dodecahedron manifold (see also [10]).

## 6. Construction of the hyperbolic football manifold $M$ = $\mathbf{H}^3$/ G

*Now we shorten our discussion, because of the strict analogy of the previous example in Sect.4.* The half-turn with axis $D_1D_2$ with an orthogonal plane to $A_0A_3$ halves the simplex (5, 3, 5). If we reflect this half simplex around vertex $A_3$, then the spherical icosahedral or dodecahedral group [2, 3, 5] (by its 120 elements) just derives the football polyhedron {5, 6, 6}. Moreover, the half-turn extended reflection group **N**:= **2 ° L** (extended *Lanner group* in the diagram of Fig.5) fills up the space $\mathbf{H}^3$ with the congruent copies of our football.

Now we define a new fixed-point-free group (denoted again by) **G** by pairing the faces of our football $\tilde{F}$ = {5, 6, 6} as fundamental domain of **G**. Again, three edges have to fall in each equivalence class, since 3 edge domains meet at every edge of the space tiling with joins pentagon-hexagon, hexagon-hexagon, hexagon-pentagon.

We start with three equivalent directed edges denoted by ⟶ . We choose the first one incident with the face assigned by $a^{-1}$. Then we introduce the screw motion **a** mapping the first edge ⟶ into the second one, the $a^{-1}$-hexagon onto the $a$-hexagon and $F_{\mathbf{G}}$ onto its **a**-image $F_{\mathbf{G}}^{\mathbf{a}}$ along the $a$-face of $F_{\mathbf{G}}$. We define the screw motion **b**, mapping the second edge ⟶ onto the third one, the $b^{-1}$ pentagon onto the $b$-pentagon and $F_{\mathbf{G}}$ onto its **b**-image $F_{\mathbf{G}}^{\mathbf{b}}$ adjacent to $F_{\mathbf{G}}$ along its $b$-face. Now, we are required to associate the hexagon $b^{-1}\,a^{-1} = (ab)^{-1}$ with the hexagon $ab$ by the product map **ab**. Indeed, the first edge ⟶ is surrounded in the $F_{\mathbf{G}}$-tiling as follows first: $F_{\mathbf{G}}$;

second: along its $a^{-1}$-hexagon comes $F_{\mathbf{G}}^{\mathbf{a}^{-1}}$ (the $\mathbf{a}^{-1}$-image of $F_{\mathbf{G}}$);

third: along the $\mathbf{a}^{-1}$-image of $b^{-1}$-pentagon comes $F_{\mathbf{G}}^{\mathbf{b}^{-1}\mathbf{a}^{-1}}$ ;

finally: along the $\mathbf{b}^{-1}\mathbf{a}^{-1}$-image of $ab$-pentagon comes $F_{\mathbf{G}}^{(\mathbf{ab})\mathbf{b}^{-1}\mathbf{a}^{-1}} = F_{\mathbf{G}}$ which is the starting polyhedron. At the same time we have explained our notations and the fixed-point-free action of **G** at inner points of the starting edge ⟶ and its images.

Now, the process is straightforward. We take the common edge of the hexagons $a^{-1}$ and $a$. Then the $\mathbf{a}^{-1}$-image and **a**-image have been determined and the directed edge class 1 is defined. Moreover, the $a^{-2}$-pentagon and $a^2$-pentagon are assigned and the screw motion $\mathbf{a}^2$ can be given as a new identifying mapping. So, we have already guaranteed the free action of **G** at inner points of the 1-edges, etc. We pick out a new directed edge in common with faces already paired, determine its edge class, either define a new identifying mapping with its paired faces or get a relation. A relation is either trivial (e.g., at the edge classes 3 and 5), or non-trivial (e.g., at 23 and 26), or a consequence of previous relations (e.g., at 27-29). The process in Fig. 5 is briefly written down as follows:

$$\mathbf{-\!\!\blacktriangleright\!\!-: a,\ b,\ ab,} \quad 1\!: \mathbf{a^2}, \quad 2\!: \mathbf{a^3} \quad 3\!: \mathbf{a \cdot a^2 \cdot a^{-3} = 1-}, \quad 4\!: \mathbf{a^3b}, \quad 5\!: -, \quad 6\!: \mathbf{ab^2},$$
$$7\!: \mathbf{a^3b^{-1}}, \quad 8\!: \mathbf{a^3b^{-1}a^{-1} =: c}, \quad 9\!: -, \quad 10\!: -, \quad 11\!: \mathbf{ba^{-2}b^2} \quad 12\!: \mathbf{a^3b^{-1}a^2},$$
$$13\!: -, \quad 14\!: \mathbf{a^3b^{-1}a^2b^{-1} =: d}, \quad 15\!: -, \quad 16\!: -, \quad 17\!: \mathbf{a^3b^{-1}a^3b =: e}, \quad 18\!: -,$$
$$19\!: \mathbf{a^3b^{-1}a^2b^{-2}a^{-1} =: f}, \quad 20\!: -, \quad 21\!: \mathbf{b^{-2}a^2b^{-1}a^3b}, \quad 22\!: -,$$
$$23\!: \mathbf{(a^3b^{-1}a^2b^{-2}a^{-1})(b^{-2}a^{-1})(b^{-2}a^2b^{-1}a^3b) = 1},$$
$$24\!: \mathbf{ba^{-2}ba^{-3}ba^{-2}b^2 =: g},$$
$$25\!: \mathbf{(ba^{-2}ba^{-3}ba^{-2}b^2)(b^{-2}a^2b^{-1}a^3b)(b^{-2}a^2b^{-1}) = 1}$$
$$26\!: \mathbf{(a^3b^{-1}a^3b)(ba^{-2}ba^{-3})(ba^{-2}ba^{-3}ba^{-2}b^2) = 1},$$
$$27\!: \mathbf{(a^3b^{-1}a^2)(ab^2a^{-2}ba^{-3})(ba^{-2}ba^{-3}ba^{-2}b^2) \overset{26}{=} 1-},$$
$$28\!: \mathbf{(a^3b^{-1}a^3b)(ba^{-2}ba^{-3}ba^{-2}b^2)(b^{-1}a^{-3}ba^{-2}b^2) \overset{26}{=} 1-},$$
$$29\!: \mathbf{(a^3b^{-1}a^2b^{-2}a^{-1})(b^{-2}a^2b^{-1}a^3b)(ba^{-2}ba^{-3}ba^{-2}b^2) \overset{23}{=} (b^{-1}a^{-3}ba^{-2}b^2a)(a^2b^{-1}a^3b)(ba^{-2}ba^{-3}ba^{-2}b^2) \overset{26}{=} 1-}, \tag{6.1}$$

In fact, the identifications satisfy the requirements for edge classes. It is already obvious, but we also realize that four vertices of $F_G$ are contained in each class of **G** equivalence, four edge classes start or end with them. To summarize, we formulate

**Theorem 1**. *The Archimedean solid (5, 6, 6) can be equipped with face identifications and locally hyperbolic metric so that it becomes an orientable compact hyperbolic space form M* $:= \mathbf{H}^3/\mathbf{G} \sim \tilde{F}_G := (5, 6, 6)$*. The fundamental group* **G** *is generated by the screw motions* **a** *and* **b** *(in Fig. 5) and the relations of* **G** *are given at edge classes 23 and 26. The first homology group of M is* $Z_{14}$*, by* **G**/[**G**, **G**]*, i.e., taking* **G** *to be commutative.*

**Proof**: The first part is proved by the constructions and the so-called Poincaré theorem as mentioned. The relations to edge class 23 and to 26 provide us just the presentation of **G**. The first homology group is

$$H_1(M) = \mathbf{G}/[\mathbf{G}, \mathbf{G}] = Z_{14} \tag{6.2}$$

the cyclic group mod 14. It can be easily calculated from the presentation by taking the commutator factor group of **G**. From the relations 23 of (6.1) we get $a^8b^{-7} = 1$, from 26 of (6.1) $a^6\, b^7 = 1$, by so-called abelianization. That means $a^2 = 1$, $b^7 = 1$. Q. e. d.

Now, let us consider the super group **N** of **G** which can be written as a decomposition

$$\mathbf{N} = \mathbf{C}_3 * \mathbf{G} \tag{6.3}$$

by (6.4). That means, each element **n** of **N** can be uniquely written in the form $\mathbf{n} = \mathbf{c}_3 \cdot \mathbf{g}$ with $\mathbf{c}_3 \in \mathbf{C}_3$, $\mathbf{g} \in \mathbf{G}$. This makes possible to express the generators of **G** by the generators of **N** (see Fig. 5)

$$\mathbf{a}^{-1} = \mathbf{r\, m_0\, m_1\, m_2\, m_1}, \quad \mathbf{b = m_3\, m_0\, m_2\, m_1\, m_0\, m_1}\,. \tag{6.4}$$

Imagine the vertex $A_2$ of $F_G$ in the middle of the face $a$ and $D_1$ at the end of edge $-\!\!\blacktriangleright\!\!-$.

Of course, the expressions in (6.4) are not unique. We could determine the distance and angle of our generating screw motions **a, b**.

*Remarks*: We have some freedom in halving $F_L$ to $F_N$. So, we could obtain other (also non-convex) fundamental domains for the group **G**.

1. In [8] we computed important data for this football manifold. It has seemingly extremal inscribed ball and extremal circumscribed ball. Its maximal ball packing density is 0.77147... and minimal ball covering density is 1.36893 …, respectively (the corresponding ball volume is related to the polyhedron volume $\mathrm{Vol}(\tilde{F})$. Both are conjectured (by the authors) best possible among all ball packings and coverings, respectively, of hyperbolic space $\mathbf{H}^3$. These extrema may involve applications in the "experimental" crystallography as well. Fullerens as "idealized constructions" initiated by Buckminster Fuller, an architect and discoverer, without any knowledge on hyperbolic geometry, have got surprising realizations [5].

2. We observe that the whole Lanner simplex (5, 3, 5), in Fig. 5 under reflections around $A_3$ with 120 group elements, generate a regular hyperbolic dodecahedron with dihedral angles $2\pi/5$. That fills $\mathbf{H}^3$ with its congruent copies, 5 dodecahedra meet at every edge. Imagine that we can construct a dodecahedron space form $\tilde{F}_\mathbf{D} = \mathbf{H}^3/\mathbf{D} = \boldsymbol{M}_1$ covering two-times our former football manifold, again with fixed-point-free orientable face pairing group **D**. The result is seen in Fig. 6 with the face paired dodecahedron $\tilde{F}_\mathbf{D}$ as fundamental domain. Observe that five edges fall in every **D**-equivalence class, indeed by the dihedral angle $2\pi/5$ above.

   The second author found in [10] all the 12 analogous hyperbolic dodecahedron face pairings (up to dodecahedron symmetries) with his systematic computer program.

## 7. Hyperbolic Cobweb Manifolds, as infinite series of possible material structures?!

Our next section follows the former ideas in Figures 7-8. This is a new topic by our papers [7] and [8]. We indicate how to construct an infinite series of orientable compact hyperbolic manifolds (space forms).

The characteristic simplex (orthoscheme) of the regular dodecahedron (5, 3, 5) and its Coxeter-Schläfli matrix in formula (5.1) can be extended, where the vertices $A_3$ and $A_0$ are outer ideal points of the hyperbolic space $\mathbf{H}^3$. This can be achieved, if we take new natural parameters $(u, v, w)$, above, $3 \le u,v,w \in \mathbf{N}$ (natural numbers), so that

$$\pi/u + \pi/v < \pi/2 \quad \text{and} \quad \pi/v + \pi/w < \pi/2, \tag{7.1}$$

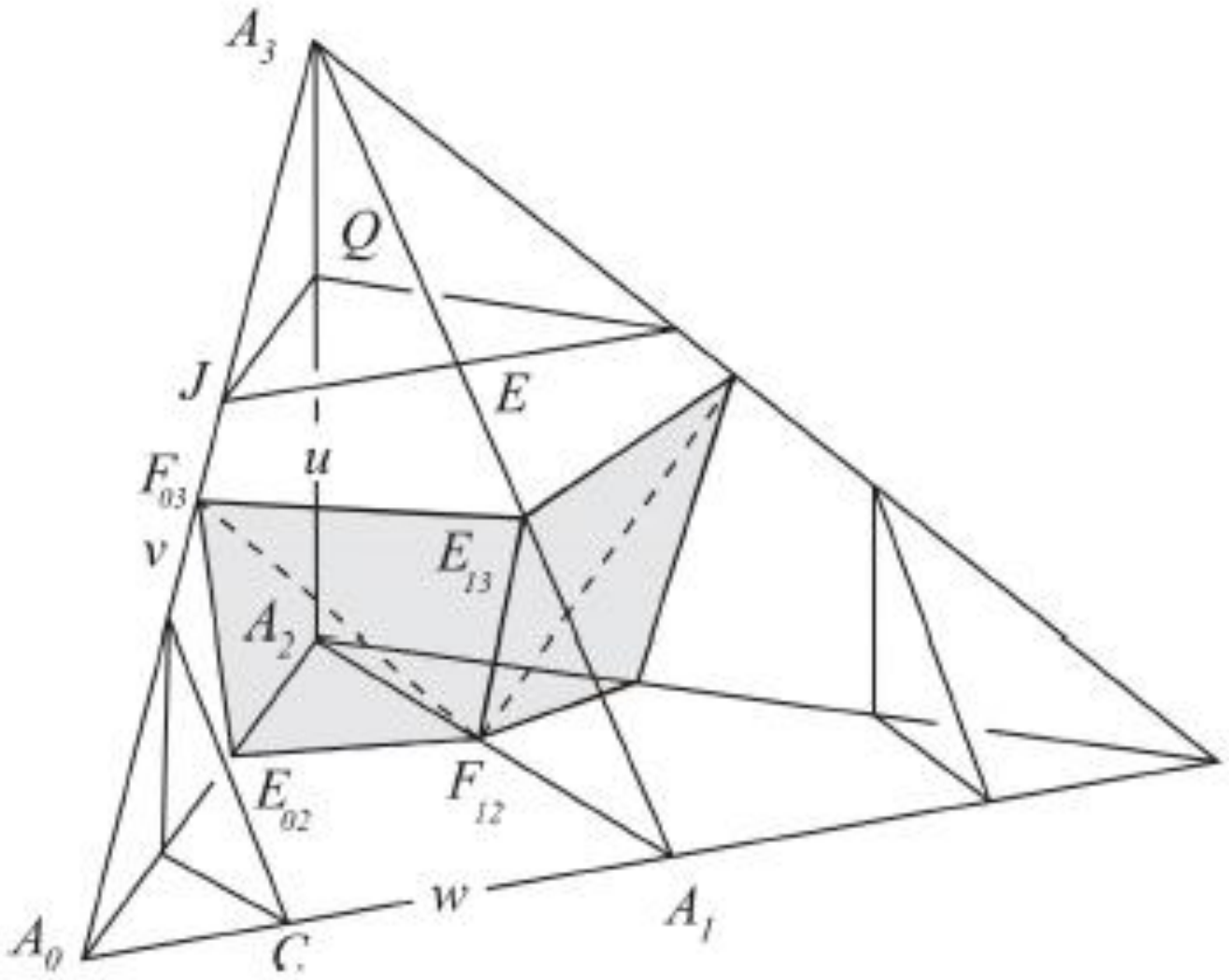

Figure 7: Fundamental domains for the half orthoscheme $W_{uvw} = W_{666}$ and for the gluing procedure at point $Q$ for getting cobweb manifold $Cw(6; 6; 6)$.

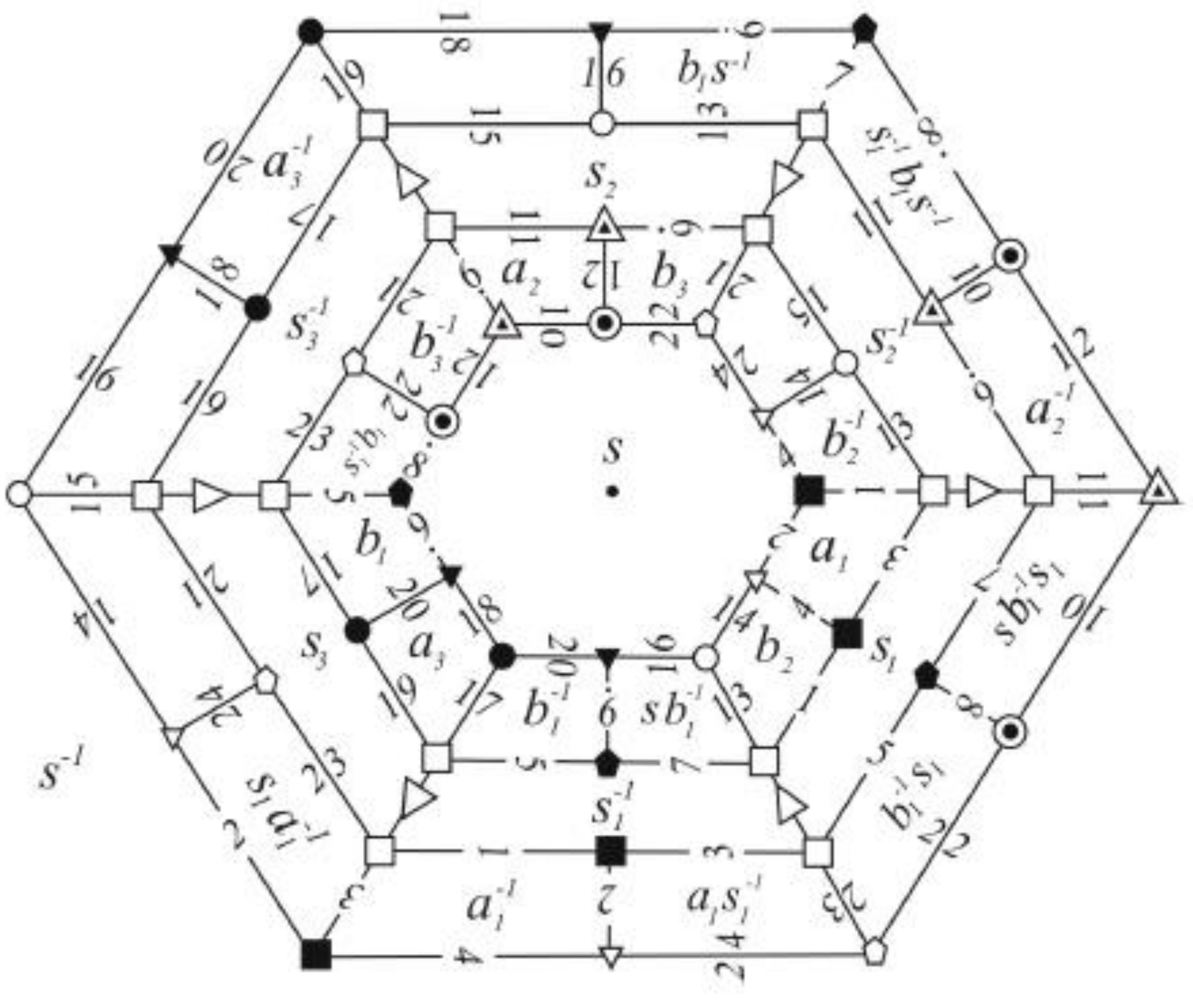

Figure 8: The hyperbolic cobweb manifold $Cw(6; 6; 6)$ with complete (symbolic) face pairing.

Imagine (Fig. 7) these dihedral face angles as $\pi/u = \beta^{01}$, $\pi/v = \beta^{12}$, $\pi/2 = \beta^{02}$ at vertex $A_3$, then the subdeterminant of (5.1) (by rows and columns 0, 1, 2) will be negative; and similarly $\pi/v = \beta^{12}$, $\pi/w = \beta^{23}$, $\pi/2 = \beta^{13}$ at vertex $A_0$. The complete determinant of (5.1) remains also negative:

$$B = \det(b^{ij}) = \sin^2\frac{\pi}{u}\sin^2\frac{\pi}{w} - \cos^2\frac{\pi}{v} < 0, \text{ i.e. } \sin\frac{\pi}{u}\sin\frac{\pi}{w} - \cos\frac{\pi}{v} < 0. \tag{7.2}$$

Then we cut the simplex by their polar planes $a_3$ and $a_0$, respectively, to get a truncated orthoscheme as generalization. The complete truncated orthoscheme group will be generated by these 5 plane reflections $\boldsymbol{b}^0$, $\boldsymbol{b}^1$, $\boldsymbol{b}^2$, $\boldsymbol{b}^3$, $\boldsymbol{a}_3$, $\boldsymbol{a}_0$, with a compact fundamental domain $W_{uvw}$ filling the hyperbolic space $\mathbf{H}^3$. In Fig. 7 we further specialize the situation by taking

$$u = v = w = 6, \text{ or more generally}$$

$$u = v = w = 2\cdot z, \text{ where } 3 \le z \text{ is odd natural number.} \tag{7.3}$$

Then a further half-turn with axis $F_{03}F_{12}$ and a halving plane $F_{03}E_{02}F_{12}E_{13}$ (with great freedom) can be introduced (analogously as at the football manifold construction).

Now, reflect this half $W_{666}$ domain around point $Q = A_2A_3 \cup a_3$. Then we obtain a cobweb solid $Cw(6, 6, 6)$, as in Fig.8, with two basis faces, 6 (or $2\cdot z$) hexagons in the middle, furthermore 24 (or $8\cdot z$ deltoids at the two basis faces.

Again, we shall introduce face pairing isometries, so that we get a fixed point free group **Cw** and a compact hyperbolic manifolds $Cw = \mathbf{H}^3/\mathbf{Cw}$. The result is completely seen in Fig. 8.

We emphasize only the most important arguments (by Figures 7-8).

We started by the 3 (= $z$) half screws

$$\mathbf{s}_1\colon s_1^{-1} \to s_1, \quad \mathbf{s}_2\colon s_2^{-1} \to s_2, \quad \mathbf{s}_3\colon s_3^{-1} \to s_3, \tag{7.4}$$

in the middle of $Cw$ around in Fig.8, so that the 6 edges with arrow —▷—, with dihedral angles $\pi/3$ fall into one equivalence class.

Then we concentrate on the odd numbered edge classes 1 – 23 (corresponding to $F_{03}E_{02} = F_{03}E_{13}$ in Fig 7) each containing 3 edges; moreover on the even numbered edge classes 2 – 24 (corresponding to $F_{12}E_{02} = F_{12}E_{13}$ in Fig 7) each containing 3 edges, too. This is because 3 copies of $Cw$ have to meet at these image edges in the space filling, respectively. To this the deltoid faces should be chosen adequately. First to

$$1\colon \mathbf{a}_1\colon a_1^{-1} \to a_1,\ \mathbf{b}_2\colon b_2^{-1} \to b_2 \text{ on faces } s_1^{-1} \text{ and } s_1, \text{ then to}$$

$$2\colon \mathbf{s}\colon s^{-1} \to s, \text{ as basis face pairing} \tag{7.5}$$

will be introduced. This becomes so luckily that the procedure smoothly goes to the end without contradiction.

Observe also the equivalence classes of vertices. At each vertex of $Cw(6,6,6)$ we shall have a ball-like neighbourhood, if we compare the situation with Fig. 7.

**Theorem 2.** *The hyperbolic Cobweb Manifold Cw =* $\mathbf{H}^3$**/Cw** *has been constructed.*

The cyclical symmetry (and our "fortune") will guarantee that this procedure can be extended for 2·$z$ "cobweb prisms", $3 \leq z$ is odd number. This will be a next publication on the base of [7] and [8], indicated here as well.

It is an open problem, what happens for 4$p$-side "cobweb prisms" ($2 \leq p \in \mathbf{N}$)?

*And, is there any application in "experimental crystallography"?*

# Az {5, 6, 6} futball és geometriája: a játékszertől a fullerénekig és tovább

Molnár Emil – Prok István – Szirmai Jenő
Budapesti Műszaki és Gazdaságtudományi Egyetem, Matematika Intézet, Geometria Tanszék

*Összefoglaló:*

Ez a "tudomány-népszerűsítő" előadás a szabályos sokszögekből indul ki. Aztán Prok I. honlapját bemutatva, a Platón-féle szabályos és az Archimédesz-féle félig szabályos testekkel, ezek szimmetriáival folytatja. A gömbi geometria síktükrözései és az úgy-nevezett alaptartományok elemzése lesz a fő eszközünk. Például a kocka középpontjában $\pi/4 = 45°$, $\pi/3 = 60°$ és $\pi/2 = 90°$ lapszögű szimmetriasíkok találkoznak. Ezekre a síkokra tükrözve, ennek a karakterisztikus gömbháromszögnek 48 példánya kövezi ki a gömböt, amit úgy is mondunk: a kocka szimmetriacsoportja 48 elemű és 3 síktükrözés generálja. A többi szabályos és félig szabályos testet, sőt síkbeli szabályos mintákat is jellemezhetünk így. A karakterisztikus háromszög szögösszegének a $\pi = 180°$-tól való eltérése jellemzi a szabályos testeket, ha az eltérés pozitív; a kockánál pl. $\pi/3 + \pi/4 + \pi/2 - \pi = \pi/12$. Az euklideszi mintáknál ez a szögkülönbség 0. De elképzelhető negatív eltérés is, ha a háromszög szögösszege kisebb 180°-nál. Ez jellemzi a Bolyai-Lobacsevszkij-féle hiperbolikus sík mintáit, kövezéseit. A tér kövezései esetében is fontos szerepet játszanak a tükrözések, melyeket a lineáris algebra eszközeivel terjeszthetünk ki L. Schläfli és H.S.M. Coxeter nyomán.

Kiderül, hogy az {5, 6, 6} szimbólumú futball-labda, melynél minden csúcsban egy szabályos (gömbi) ötszög és két szabályos (gömbi) hatszög találkozik, poliéder – azaz síklapú test – formájában, egybevágó példányaival nem tudja kitölteni euklideszi terünket, de a Bolyai Lobacsevszkij-féle hiperbolikus teret ki tudja tölteni (ahogy ezt Molnár E. 1988-ban észrevette és publikálta. Sőt a kitöltés fixpont-mentes egybevágóságokkal történhet: A futball-labdák olyanok, mintha egy véges térben lennénk, ahol minden pontnak hiperbolikus labda környezete van. Úgy tűnik, hogy az anyag tudományokban előtérbe került *fullerének*, azaz $C_{60}$ molekulák ezt a nagyon „sűrű hiperbolikus anyagot" követik (melynek további szélső-érték tulajdonságai Szirmai J. és Molnár E. kutatásaiban is előjönnek).

Tehát Bolyai János és N.I. Lobacsevszkij hiperbolikus geometriájának kristálytani alkalmazásai is lehetnek az eddigi csillagméretű vonatkozások mellett.